\documentclass[10pt]{amsart}
\usepackage{amssymb}
\theoremstyle{definition}
\newtheorem{definition}{Definition}[section]

\newtheorem{remark}[definition]{Remark}
\theoremstyle{plain}
\newtheorem{theorem}[definition]{Theorem}
\newtheorem{lemma}[definition]{Lemma}

\newtheorem{assumption}{Assumption}

\def\m#1{\theta_{#1}}
\def\M{\Theta}
\def\oo{\overline}

\title[A class of Bol loops]{A class of Bol loops with a subgroup of index two}

\author[Petr~Vojt\v{e}chovsk\'y]{Petr Vojt\v{e}chovsk\'y}
\email{petr@math.du.edu}
\address{Department of Mathematics, University of Denver, 2360 S Gaylord St,
Denver, CO 80208, USA}

\begin{document}

\begin{abstract}Let $G$ be a finite group and $C_2$ the cyclic group of order
$2$. Consider the $8$ multiplicative operations $(x,y)\mapsto (x^iy^j)^k$,
where $i$, $j$, $k\in\{-1,\,1\}$. Define a new multiplication on $G\times C_2$
by assigning one of the above $8$ multiplications to each quarter
$(G\times\{i\})\times(G\times\{j\})$, for $i$, $j\in C_2$. We describe all
situations in which the resulting quasigroup is a Bol loop. This paper also
corrects an error in P.~Vojt\v{e}chovsk\'y: On the uniqueness of loops
$M(G,2)$. \vskip 1mm

 \vskip 1mm

\end{abstract}

\maketitle

\section{Introduction}

\noindent Let $G$ be a finite group. Consider the two maps $\m{yx}$,
$\m{xy^{-1}}: G\times G\to G\times G$ defined by $\m{yx}(a,b)=(b,a)$,
$\m{xy^{-1}}(a,b)=(a,b^{-1})$. The group $\M=(\M,\circ)$ generated by $\m{yx}$
and $\m{xy^{-1}}$ consists of eight maps, and is isomorphic to the quaternion
group. We will denote elements of $\M$ suggestively by $\m{xy}$, $\m{xy^{-1}}$,
$\m{x^{-1}y}$, $\m{x^{-1}y^{-1}}$, $\m{yx}$, $\m{yx^{-1}}$, $\m{y^{-1}x}$ and
$\m{y^{-1}x^{-1}}$. For instance, $\m{y^{-1}x}$ is the map defined by
$\m{y^{-1}x}(a,b)=(b^{-1},a)$.

Let us identify $\m{uv}\in \M$ with $\Delta\m{uv}$, where $\Delta:G\times G\to
G$ is given by $\Delta(a,b)=ab$. Thanks to this double perspective, each
$\m{uv}\in \M$ determines a multiplication on $G$, yet it is possible to
compose the multiplications.

Let $\oo{G}=\{\oo{g};\;g\in G\}$ be a copy of $G$. Given four maps $\alpha$,
$\beta$, $\gamma$, $\delta\in \M$, define multiplication $*$ on $G\cup \oo G$
by
\begin{displaymath}
    g*h=\alpha(g,h),\quad g*\oo{h}=\oo{\beta(g,h)},\quad
    \oo{g}*h=\oo{\gamma(g,h)},\quad
    \oo{g}*\oo{h}=\delta(g,h),
\end{displaymath}
or, more precisely,
\begin{displaymath}
    g*h=\Delta\alpha(g,h),\quad g*\oo{h}=\oo{\Delta\beta(g,h)},\quad
    \oo{g}*h=\oo{\Delta\gamma(g,h)},\quad
    \oo{g}*\oo{h}=\Delta\delta(g,h),
\end{displaymath}
where $g$, $h\in G$. Note that $G*G=G=\oo{G}*\oo{G}$,
$G*\oo{G}=\oo{G}=\oo{G}*G$. Hence $G(\alpha,\beta,\gamma,\delta)=(G\cup
\oo{G},*)$ is a quasigroup with normal subgroup $G$.

Chein proved in \cite{Chein} that
\begin{displaymath}
    M(G,2)=G(\m{xy},\m{yx},\m{xy^{-1}},\m{y^{-1}x})
\end{displaymath}
is a Moufang loop for every group $G$, and that is is associative if and only
if $G$ is abelian. Many small Moufang loops are of this kind, e.g. the smallest
nonassociative Moufang loop $M(S_3,2)$ (see \cite{CP}, \cite{GMR},
\cite{Vojtechovsky}).

It is natural to ask if there are other constructions
$G(\alpha,\beta,\gamma,\delta)$ besides $M(G,2)$ that produce Moufang loops. I
have obtained the following result in \cite{PV}:

\begin{theorem}\label{Th:Moufang}
Let $G$ be a nonabelian group, and let $\alpha$, $\beta$, $\gamma$,
$\delta\in\M$ be as above. Then $G(\alpha,\beta,\gamma,\delta)$ is a Moufang
loop if and only if it is among
\begin{displaymath}
    \begin{array}{ll}
        G(\m{xy},\m{xy},\m{xy},\m{xy}), &
        G(\m{yx},\m{yx},\m{yx},\m{yx}),\\
        G(\m{xy},\m{yx^{-1}},\m{y^{-1}x},\m{x^{-1}y^{-1}}), &
        G(\m{yx},\m{x^{-1}y},\m{xy^{-1}},\m{y^{-1}x^{-1}}),\\
        G(\m{xy},\m{yx},\m{xy^{-1}},\m{y^{-1}x}), &
        G(\m{yx},\m{yx^{-1}},\m{xy},\m{x^{-1}y}), \\
        G(\m{xy},\m{x^{-1}y},\m{yx},\m{yx^{-1}}), &
        G(\m{yx},\m{xy},\m{y^{-1}x},\m{xy^{-1}}).
    \end{array}
\end{displaymath}
The four loops in the first two rows are associative and isomorphic to the
direct product of $G$ with the $2$-element cyclic group. The remaining four
loops are not associative and are isomorphic to the loop $M(G,2)$.
\end{theorem}

\section{The Mistake}

\noindent Hence Chein's construction $M(G,2)$ is the ``unique'' construction
$G(\alpha,\beta,\gamma,\delta)$ that produces nonassociative Moufang loops. I
claim in \cite{PV} that $M(G,2)$ is the ``unique'' construction
$G(\alpha,\beta,\gamma,\delta)$ that produces nonassociative Bol loops, too.
Unfortunately, this is not correct, as we shall see.

The mistake (pointed out to me by Michael Kinyon) is in Lemma 2 \cite{PV},
where I claim that any loop of the form $G(\alpha,\beta,\gamma,\delta)$ is an
inverse property loop. But I only prove in Lemma 2 \cite{PV} that such a loop
has two-sided inverses. Hence the conclusion that any Bol loop
$G(\alpha,\beta,\gamma,\delta)$ is automatically Moufang was not justified in
\cite{PV}, and is in fact not true.

This note can be considered an erratum for \cite{PV}. However, in dealing with
the Bol case, I had to develop additional techniques not found in \cite{PV}.
Moreover, several constructions obtained in this way appear to be new, and
should be of interest in the classification of small Bol loops.

\section{Reductions and Assumptions}

\noindent Our goal is to describe all Bol loops of the form
$G(\alpha,\beta,\gamma,\delta)$.

Note that $|\Delta \M|=1$ when $G$ is an elementary abelian $2$-group or when
$|G|=1$. In such a case, any $G(\alpha,\beta,\gamma,\delta)$ is equal to
$G(\m{xy},\m{xy},\m{xy},\m{xy})$, and is therefore isomorphic to $G\times C_2$.
We thus lose nothing by making this assumption:

\begin{assumption}\label{As:1} $|G|>1$ and $G$ is not an
elementary abelian $2$-group.
\end{assumption}

Although $G(\alpha,\beta,\gamma,\delta)$ does not have to be a loop, it is not
hard to determine when it is (cf. Lemma 1 \cite{PV}):

\begin{lemma}\label{Lm:Loop}
The quasigroup $(M,*)=G(\alpha,\beta,\gamma,\delta)$ is a loop if and only if
\begin{equation}\label{Eq:Loop}
    \alpha\in\{\m{xy},\m{yx}\},\quad
    \beta\in\{\m{xy},\m{yx},\m{yx^{-1}},\m{x^{-1}y}\},\quad
    \gamma\in\{\m{xy},\m{yx},\m{y^{-1}x},\m{xy^{-1}}\}.
\end{equation}
When $(M,*)$ is a loop, its neutral element coincides with the neutral element
of $G$.
\end{lemma}
\begin{proof}
We first show that if $(M,*)$ is a loop, its neutral element $e$ coincides with
the neutral element $1$ of $G$. We have $1*1=\alpha(1,1)=1=e*1$, no matter what
$\alpha\in\M$ is. Since $(M,*)$ is a quasigroup, $1=e$ follows.

The equation $y=1*y$ holds for every $y\in G$ if and only if $y=\alpha(1,y)$
holds for every $y\in G$, which happens if and only if $\alpha$ does not invert
its second argument. (We will use this trick many times. Note how Assumption
\ref{As:1} is used.) Thus $y=1*y$ holds for every $y\in G$ if and only if
$\alpha\in\{\m{xy},\m{x^{-1}y},\m{yx},\m{yx^{-1}}\}$.

Similarly, the equation $y=y*1$ holds for every $y\in G$ if and only if
$\alpha\in\{\m{xy},\,\m{xy^{-1}},\m{yx},\m{y^{-1}x}\}$. Altogether, $y=y*1=1*y$
holds for every $y\in G$ if and only if $\alpha\in\{\m{xy},\m{yx}\}$.

Following a similar strategy, $\oo{y}=1*\oo{y}$ holds for every $y\in G$ if and
only if $\beta\in\{\m{xy},\m{yx},\m{yx^{-1}},\m{x^{-1}y}\}$, and
$\oo{y}=\oo{y}*1$ holds for every $y\in G$ if and only if $
\gamma\in\{\m{xy},\m{yx},\m{y^{-1}x},\m{xy^{-1}}\}$.
\end{proof}

Since we are only interested in loops here, we assume:

\begin{assumption}\label{As:2}
The maps $\alpha$, $\beta$, $\gamma$ are as in $(\ref{Eq:Loop})$.
\end{assumption}

For a groupoid $(A,\circ)$, let $A^{\mathrm op}=(A,\circ^{\mathrm op})$ be the
\emph{opposite} of $A$, defined by $x\circ^{\mathrm op} y = y\circ x$. Then
\begin{eqnarray}
    G(\alpha,\beta,\gamma,\delta)^{\mathrm op} &=&
    G(\m{yx}\alpha,\m{yx}\gamma,\m{yx}\beta,\m{yx}\delta),\label{Eq:Op1}\\
    G(\alpha,\beta,\gamma,\delta) &=& G^{\mathrm op}
    (\m{yx}\alpha,\m{yx}\beta,\m{yx}\gamma,\m{yx}\delta).\label{Eq:Op2}
\end{eqnarray}
It is not necessarily true that $G(\alpha,\beta,\gamma,\delta)^{\mathrm op}$ is
isomorphic to $G(\alpha,\beta,\gamma,\delta)$. However, by $(\ref{Eq:Op2})$,
any loop $G(\m{yx},\beta,\gamma,\delta)$ can also be obtained as $G^{\mathrm
op}(\m{xy},\beta',\gamma',\delta')$, for some $\beta'$, $\gamma'$, $\delta'$.
As $G$ is isomorphic to $G^{\mathrm op}$ (via $x\mapsto x^{-1})$, we postulate:

\begin{assumption}\label{As:3}
$\alpha=\m{xy}$
\end{assumption}

Our last reduction concerns the maps $\beta$ and $\gamma$.

\begin{lemma}\label{Lm:Iso} The loop $G(\alpha,\beta,\gamma,\delta)$ is
isomorphic to $G(\alpha,\beta',\gamma',\m{x^{-1}y^{-1}}\delta)$ if
\begin{eqnarray*}
    (\beta,\beta')&\in&\{(\m{xy},\m{yx^{-1}}),(\m{yx},\m{x^{-1}y}),
        (\m{x^{-1}y},\m{yx}), (\m{yx^{-1}},\m{xy})\},\\
    (\gamma,\gamma')&\in&\{(\m{xy},\m{y^{-1}x}),(\m{yx},\m{xy^{-1}}),
        (\m{xy^{-1}},\m{yx}),(\m{y^{-1}x},\m{xy})\}.
\end{eqnarray*}
\end{lemma}
\begin{proof}
Let $(M,*)=G(\alpha,\beta,\gamma,\delta)$, and
$(M,\circ)=G(\alpha,\beta',\gamma',\m{x^{-1}y^{-1}}\delta)$. Consider the
permutation $f:M=G\cup \oo G\to M$ defined by $f(x)=x$, $f(\oo x)=\oo{x^{-1}}$,
for $x\in G$.

We show that $f$ is an isomorphism of $(M,*)$ onto $(M,\circ)$ if (and only if)
\begin{equation}\label{Eq:Aux1}
    (\Delta\beta(x,y))^{-1}=\Delta\beta'(x,y^{-1}),\quad
    (\Delta\gamma(x,y))^{-1}=\Delta\gamma'(x^{-1},y).
\end{equation}
Once we establish this fact, the proof is finished by checking that the pairs
$(\beta,\beta')$, $(\gamma,\gamma')$ in the statement of the Lemma satisfy
$(\ref{Eq:Aux1})$.

In the following computation, we emphasize by $\Delta$ the multiplication in
$M$. Let $x$, $y\in G$. Then
\begin{eqnarray*}
    f(x*y)&=&f(\Delta\alpha(x,y)) = \Delta\alpha(x,y),\\
    f(x*\oo y)&=&f(\oo{\Delta\beta(x,y)}) = \oo{(\Delta\beta(x,y))^{-1}},\\
    f(\oo x*y)&=&f(\oo{\Delta\gamma(x,y)}) = \oo{(\Delta\gamma(x,y))^{-1}},\\
    f(\oo x*\oo y)&=&f(\Delta\delta(x,y))=\Delta\delta(x,y),
\end{eqnarray*}
while
\begin{eqnarray*}
    f(x)\circ f(y)&=& x\circ y = \Delta\alpha(x,y),\\
    f(x)\circ f(\oo y)&=& x\circ \oo{y^{-1}} = \oo{\Delta\beta'(x,y^{-1})},\\
    f(\oo x)\circ f(y)&=& \oo{x^{-1}}\circ y = \oo{\Delta\gamma'(x^{-1},y)},\\
    f(\oo x)\circ f(\oo y)&=& \oo{x^{-1}}\circ\oo{y^{-1}} =
        \Delta\m{x^{-1}y^{-1}}\delta(x^{-1},y^{-1}).
\end{eqnarray*}
We see that $f(x*y)=f(x)\circ f(y)$, $f(\oo x*\oo y)=f(\oo x)\circ f(\oo y)$
always hold, and that $f(x*\oo y)=f(x)\circ f(\oo y)$, $f(\oo x*y)=f(\oo
x)\circ f(y)$ hold if $(\beta,\beta')$, $(\gamma,\gamma')$ satisfy
$(\ref{Eq:Aux1})$.
\end{proof}

Observe that for any admissible value of $\beta$, $\gamma$ and $\delta$, Lemma
\ref{Lm:Iso} provides an isomorphism of $G(\alpha,\beta,\gamma,\delta)$ that
leaves $\alpha$ intact. Furthermore, if $\gamma=\m{xy^{-1}}$, the corresponding
$\gamma'$ is equal to $\m{yx}$, and if $\gamma=\m{y^{-1}x}$, we have
$\gamma'=\m{xy}$. We can therefore assume:

\begin{assumption}\label{As:4}
    $\gamma\in\{\m{xy},\m{yx}\}$
\end{assumption}

We could have reduced $\beta$ instead of $\gamma$, but the choice we have made
will be more convenient later on.

\section{The Technique}

\noindent Suppose we want to check if $(M,*)=G(\alpha,\beta,\gamma,\delta)$
satisfies a given identity. Since the product $a*b$ in $M$ depends on whether
the elements $a$, $b$ belong to $G$ or $\oo G$, it is natural to treat the
cases separately. For instance, the \emph{left alternative law}
$a*(a*b)=(a*a)*b$ for $M$ leads to four identities $x*(x*y)=(x*x)*y$, $x*(x*\oo
y) = (x*x)*\oo y$, $\oo x*(\oo x * y) = (\oo x * \oo x) * y$, and $\oo x * (\oo
x * \oo y) = (\oo x*\oo x)*\oo y$, where $x$, $y\in G$. In turn, each of these
identities can be rewritten as an identity for $G$, using the maps $\alpha$,
$\beta$, $\gamma$, $\delta$. Here are the four left alternative identities
together with their translations:
\begin{eqnarray}
       x*(x*y)=(x*x)*y\quad &
            \alpha(x,\alpha(x,y)) = \alpha(\alpha(x,x),y),\label{Eq:LA1}\\
       x*(x*\oo y) = (x*x)*\oo y\quad &
            \beta(x,\beta(x,y)) = \beta(\alpha(x,x),y),\label{Eq:LA2}\\
        \oo x*(\oo x * y) = (\oo x * \oo x) * y\quad &
            \delta(x,\gamma(x,y)) = \alpha(\delta(x,x),y),\label{Eq:LA3}\\
        \oo x * (\oo x * \oo y) = (\oo x*\oo x)*\oo y\quad &
            \gamma(x,\delta(x,y)) = \beta(\delta(x,x),y).\label{Eq:LA4}
\end{eqnarray}
In case we need to prove that $(M,*)$ does not satisfy a given identity, it
suffices to show that any of the translated identities does not hold for $G$.
For such purposes, it is often advantageous to look at an identity that does
not involve many different maps (i.e., $(\ref{Eq:LA2})$ is preferable to
$(\ref{Eq:LA4})$, say).

More importantly, since we know nothing about $G$ besides the fact that it
satisfies Assumption \ref{As:1}, how do we decide if a given identity is true
or false in $G$? Well, if we treat the identity as an identity in a free group,
and if the identity reduces to $x=x^{-1}$, it must be false in $G$. More
identities reduce to $x=x^{-1}$ if we assume that the free group is abelian.
Instead of making this distinction for each particular identity, we will treat
the abelian and nonabelian cases separately from the start.

When $G$ is a group and $m$ a positive integer, we let
\begin{displaymath}
    G^m=\{g^m;\;g\in G\},\quad G_m=\{g\in G;\;g^m=1\}.
\end{displaymath}
Note that, in the literature, $G_m$ often denotes the set of all elements whose
order is a power of $m$. In our case, $G_m$ consists of all elements of
exponent $m$.

Because of the nature of Bol identities (see below), we will often come across
group identities involving squares. Two conditions for $G$ will then help us
characterize the groups in which such identities hold. Namely, $G^2\subseteq
Z(G)$, and $G^4=1$. The former assumption says that $G/Z(G)$ is an elementary
abelian $2$-group. The latter assumption is equivalent to $G_4=G$.

We are now ready to start the search for all Bol loops
$G(\alpha,\beta,\gamma,\delta)$. Recall that a loop is \emph{left Bol} if it
satisfies the identity
\begin{equation}\label{Eq:LB}
    x(y(xz)) = (x(yx))z,
\end{equation}
 and it is \emph{right Bol} if it satisfies the identity
$((zx)y)x = z((xy)x)$. (See \cite{Pflugfelder}.) Hence left Bol loops are
opposites of right Bol loops, and thus all right Bol loops of the form
$G(\alpha,\beta,\gamma,\delta)$ can be obtained from the left Bol ones via
$(\ref{Eq:Op1})$.

We therefore restrict our search to left Bol loops.

\section{The Abelian Case}

\noindent In this section we suppose that $G$ is a finite abelian group. Then
$|\Delta\M|=4$, and it suffices to consider multiplications $\m{xy}$,
$\m{x^{-1}y}$, $\m{xy^{-1}}$, $\m{x^{-1}y^{-1}}$.

Assumptions \ref{As:1}--\ref{As:4} are now equivalent to: \emph{$G$ is abelian,
$|G|>1$, $G$ is not an elementary abelian $2$-group, and}
\begin{displaymath}
    \alpha=\m{xy},\ \beta\in\{\m{xy},\m{x^{-1}y}\},\
    \gamma=\m{xy},\
    \delta\in\{\m{xy}, \m{x^{-1}y},\m{xy^{-1}},\m{x^{-1}y^{-1}}\}.
\end{displaymath}

Left Bol loops are left alternative, as is immediately obvious upon
substituting $y=1$ into $(\ref{Eq:LB})$. Let us therefore first describe all
left alternative loops $G(\alpha,\beta,\gamma,\delta)$.

\begin{lemma}\label{Lm:CLA}
Suppose that $G$ is an abelian group and that Assumptions
$\ref{As:1}$--$\ref{As:4}$ are satisfied. Then $G(\alpha,\beta,\gamma,\delta)$
is a left alternative loop if and only if one of the following conditions is
satisfied:
\begin{enumerate}
\item[(i)] $(\beta,\delta)\in\{(\m{xy},\m{xy}),(\m{xy},\m{x^{-1}y}),
    (\m{x^{-1}y},\m{x^{-1}y})\}$,
\item[(ii)] $G^4=1$ and $(\beta,\delta) = (\m{x^{-1}y},\m{xy})$.
\end{enumerate}
\end{lemma}
\begin{proof} Identity $(\ref{Eq:LA1})$ holds since $\alpha=\m{xy}$. Identity
$(\ref{Eq:LA2})$ clearly holds when $\beta=\m{xy}$. When $\beta=\m{x^{-1}y}$,
it becomes $x^{-1}x^{-1}y=(xx)^{-1}y$; again true.

When $\beta=\m{xy}$, identity $(\ref{Eq:LA4})$ becomes
$x\delta(x,y)=\delta(x,x)y$. Note that this identity cannot hold if $\delta$
inverts its second argument. (Since then there is $y^{-1}$ on the left, $y$ on
the right, and with $x=1$ the identity becomes $y=y^{-1}$.) On the other hand,
the identity $x\delta(x,y)=\delta(x,x)y$ holds for $\delta=\m{xy}$ and
$\delta=\m{x^{-1}y}$.

When $\beta=\m{x^{-1}y}$, $(\ref{Eq:LA4})$ becomes
$x\delta(x,y)=\delta(x,x)^{-1}y$. Again, this identity cannot hold if $\delta$
inverts its second argument. When $\delta=\m{x^{-1}y}$, it becomes $y=y$
(true). When $\delta=\m{xy}$, it becomes $xxy=(xx)^{-1}y$, which holds if and
only if $G^4=1$.

Finally, $(\ref{Eq:LA3})$ holds for $\delta\in\{\m{xy}$, $\m{x^{-1}y}\}$.
\end{proof}

\begin{theorem}\label{Th:CLB}
Assume that $G$ is an abelian group and that Assumptions
$\ref{As:1}$--$\ref{As:4}$ are satisfied. Then
$(M,*)=G(\alpha,\beta,\gamma,\delta)$ is a group if and only if it is among
\begin{eqnarray}
    &&G(\m{xy},\m{xy},\m{xy},\m{xy}),\label{Eq:CLB1}\\
    &&G(\m{xy},\m{x^{-1}y},\m{xy},\m{x^{-1}y}).\label{Eq:CLB2}
\end{eqnarray}
The group $(\ref{Eq:CLB1})$ is the direct product $G\times C_2$, and the group
$(\ref{Eq:CLB2})$ is isomorphic to the Chein loop $M(G,2)$.

Furthermore, $(M,*)$ is a nonassociative left Bol loop if and only if $(M,*)$
is a left Bol loop that is not Moufang if and only if $G^4=1$ and $(M,*)$ is
among
\begin{eqnarray}
    &&G(\m{xy},\m{xy},\m{xy},\m{x^{-1}y}),\label{Eq:CLB3}\\
    &&G(\m{xy},\m{x^{-1}y},\m{xy},\m{xy}).\label{Eq:CLB4}
\end{eqnarray}
For a given $G$, the two loops $(\ref{Eq:CLB3})$, $(\ref{Eq:CLB4})$ are not
isomorphic.
\end{theorem}
\begin{proof}
There cannot be more than $4$ left Bol loops $G(\alpha,\beta,\gamma,\delta)$
satisfying Assumptions $\ref{As:1}$--$\ref{As:4}$, since there are only $4$
such left alternative loops, by Lemma \ref{Lm:CLA}.

Clearly, $(\ref{Eq:CLB1})$ is the direct product $G\times C_2$. We claim that
the loop $(\ref{Eq:CLB2})$ is isomorphic to $M(G,2)$. To see that, write
$M(G,2)$ as $G(\m{xy},\m{xy},\m{xy^{-1}},\m{xy^{-1}})$, and apply a suitable
isomorphism of Lemma \ref{Lm:Iso} to it. Since $G$ is abelian,
$(\ref{Eq:CLB2})$ is associative by \cite{Chein}.

Consider one of the flexible identities $\oo x*(\oo y*\oo x) = (\oo x * \oo
y)*\oo x$. It translates into $\gamma(x,\delta(y,x))=\beta(\delta(x,y),x)$.
This becomes $xy^{-1}x=x^{-1}yx$ for $(\ref{Eq:CLB3})$, and $xyx=y^{-1}x^{-1}x$
for $(\ref{Eq:CLB4})$; both false (let $y=1$).

Thus neither of the loops $(\ref{Eq:CLB3})$, $(\ref{Eq:CLB4})$ is flexible, and
hence neither is a Moufang loop. We must now check that $(\ref{Eq:CLB3})$,
$(\ref{Eq:CLB4})$ are left Bol. This follows by straightforward calculation:

The left Bol identity for $(M,*)$ translates into $8$ identities for $G$. Here
they are: {\small
\begin{eqnarray}
   x*(y*(x*z)) = (x*(y*x))*z &&
   \alpha(x,\alpha(y,\alpha(x,z))) = \alpha(\alpha(x,\alpha(y,x)),z)
   \label{Eq:LB1}\\
   x*(y*(x*\oo z)) = (x*(y*x))*\oo z &&
   \beta(x,\beta(y,\beta(x,z))) = \beta(\alpha(x,\alpha(y,x)),z)
   \label{Eq:LB2}\\
   x*(\oo y*(x*z)) = (x*(\oo y*x))*z &&
   \beta(x,\gamma(y,\alpha(x,z))) = \gamma(\beta(x,\gamma(y,x)),z)
   \label{Eq:LB3}\\
   \oo x*(y*(\oo x*z)) = (\oo x*(y*\oo x))*z &&
   \delta(x,\beta(y,\gamma(x,z))) = \alpha(\delta(x,\beta(y,x)),z)
   \label{Eq:LB4}\\
   x*(\oo y*(x*\oo z)) = (x*(\oo y*x))*\oo z &&
   \alpha(x,\delta(y,\beta(x,z))) = \delta(\beta(x,\gamma(y,x)),z)
   \label{Eq:LB5}\\
   \oo x*(y*(\oo x*\oo z)) = (\oo x*(y*\oo x))*\oo z &&
   \gamma(x,\alpha(y,\delta(x,z))) = \beta(\delta(x,\beta(y,x)),z)
   \label{Eq:LB6}\\
   \oo x*(\oo y*(\oo x*z)) = (\oo x*(\oo y*\oo x))*z &&
   \gamma(x,\delta(y,\gamma(x,z))) = \gamma(\gamma(x,\delta(y,x)),z)
   \label{Eq:LB7}\\
   \oo x*(\oo y*(\oo x*\oo z)) = (\oo x*(\oo y*\oo x))*\oo z &&
   \delta(x,\gamma(y,\delta(x,z))) = \delta(\gamma(x,\delta(y,x)),z)
   \label{Eq:LB8}
\end{eqnarray}}
Since the only nontrivial multiplication in $(\ref{Eq:CLB3})$ is $\delta$, it
suffices to verify identities $(\ref{Eq:LB4})$--$(\ref{Eq:LB8})$ for
$(\ref{Eq:CLB3})$. We obtain $x^{-1}yxz=x^{-1}yxz$ (true),
$xy^{-1}xz=(xyx)^{-1}z$ (true if and only if $G^4=1$), $xyx^{-1}z=x^{-1}yxz$
(true), $xy^{-1}xz=xy^{-1}xz$ (true), and $x^{-1}yx^{-1}z=(xy^{-1}x)^{-1}z$
(true), respectively.

Since the only nontrivial multiplication in $(\ref{Eq:CLB4})$ is $\beta$, it
suffices to verify identities $(\ref{Eq:LB2})$--$(\ref{Eq:LB6})$ for
$(\ref{Eq:CLB4})$. We obtain $x^{-1}y^{-1}x^{-1}z=(xyx)^{-1}z$ (true),
$x^{-1}yxz=x^{-1}yxz$ (true), $xy^{-1}xz=xy^{-1}xz$ (true),
$xyx^{-1}z=x^{-1}yxz$ (true), $xyxz=(xy^{-1}x)^{-1}z$ (true if and only if
$G^4=1$), respectively.

It remains to check that $(\ref{Eq:CLB3})$ is not isomorphic to
$(\ref{Eq:CLB4})$. To see that, notice that all elements $\oo x$ are of order
$2$ in $(\ref{Eq:CLB3})$ (since $\delta=\m{xy^{-1}}$), while the order of $\oo
x$ is bigger or equal to the order of $x$ in $(\ref{Eq:CLB4})$ (since the $n$th
power of $\oo x$ in $(M,*)$ can be written as $\oo x*(\oo x*(\dots))$, and
$\gamma=\delta=\m{xy}$).
\end{proof}

\begin{remark} By the Fundamental theorem for finite abelian groups
\cite{Rotman}, the groups required for $(\ref{Eq:CLB3})$, $(\ref{Eq:CLB4})$ are
exactly the groups of the form $(C_4)^m\times (C_2)^n$, where $m>0$ and $n\ge
0$.
\end{remark}

\begin{remark}\label{Rm:Middle}
The smallest left Bol loops that are not Moufang are of order $8$. In fact,
there are $6$ such loops, up to isomorphism \cite{Burn}. Five of them contain a
subgroup isomorphic to $C_4$. (See \cite{Burn}, or use the LOOPS package
\cite{LOOPS} for GAP, where the six Bol loops can be obtained via
\texttt{BolLoop(8,i)}, with $1\le$\texttt{i}$\le 6$.) Our constructions
$(\ref{Eq:CLB3})$ and $(\ref{Eq:CLB4})$ yield two of these five loops. We can
of course write down their multiplication tables easily:
\begin{displaymath}
    \begin{array}{cc}
        \begin{array}{cccc|cccc}
            1&2&3&4&5&6&7&8\\
            2&3&4&1&6&7&8&5\\
            3&4&1&2&7&8&5&6\\
            4&1&2&3&8&5&6&7\\
            \hline
            5&6&7&8&1&2&3&4\\
            6&7&8&5&4&1&2&3\\
            7&8&5&6&3&4&1&2\\
            8&5&6&7&2&3&4&1
        \end{array}\quad\quad\quad
        \begin{array}{cccc|cccc}
            1&2&3&4&5&6&7&8\\
            2&3&4&1&8&5&6&7\\
            3&4&1&2&7&8&5&6\\
            4&1&2&3&6&7&8&5\\
            \hline
            5&6&7&8&1&2&3&4\\
            6&7&8&5&2&3&4&1\\
            7&8&5&6&3&4&1&2\\
            8&5&6&7&4&1&2&3
        \end{array}
    \end{array}
\end{displaymath}
\end{remark}

\begin{remark} Dr\'apal showed in \cite{Drapal} that proximity of group
multiplication tables implies proximity of algebraic properties. More
precisely, he defined the \emph{distance} $d(*,\circ)$ of two groups $(G,*)$,
$(G,\circ)$ with the same underlying set $G$ as the cardinality of $\{(g,h)\in
G\times G;\;g*h\ne g\circ h\}$, and showed that if $d(*,\circ)<|G|^2/9$ then
$(G,*)$ is isomorphic to $(G,\circ)$. If $(G,*)$ (and thus $(G,\circ)$) is a
$2$-group, the isomorphism follows already from $d(*,\circ)<|G|^2/4$ (cf.
\cite{Drapal2}). Dr\'apal and the author conjecture in \cite{AlesPetr} that the
same is true for Moufang $2$-loops, i.e., $(G,*)\cong(G,\circ)$ if $(G,*)$,
$(G,\circ)$ are Moufang $2$-loops satisfying $d(*,\circ)<|G|^2/4$.

Note that the distance of any of the two nonassociative left Bol loops in
Remark \ref{Rm:Middle} from the canonical multiplication table of the direct
product $C_4\times C_2$ is $8$, i.e., only $1/8\cdot 8^2$. Hence the conjecture
cannot be generalized from Moufang to Bol loops.
\end{remark}

\section{The Nonabelian Case}

\noindent In this section we suppose that $G$ is a finite nonabelian group.
This assumption has some consequences on the validity of identities. For
instance, while $xy^{-1}=y^{-1}x$ holds in the abelian case, it is
\emph{always} false in this section, no matter what $G$ is.

We stick to the same strategy as in the abelian case.

\begin{lemma}\label{Lm:LA}
Assume that $G$ is a nonabelian group, and that Assumptions
$\ref{As:1}$--$\ref{As:4}$ are satisfied. Then $G(\alpha,\beta,\gamma,\delta)$
is a left alternative loop if and only if one of the following conditions is
satisfied:
\begin{enumerate}
\item[(i)] $(\gamma,\delta)\in\{(\m{xy},\m{x^{-1}y}),(\m{yx},\m{yx^{-1}})\}$,

\item[(ii)] $(\beta,\gamma,\delta) =(\m{xy},\m{xy},\m{xy})$,

\item[(iii)] $G^2\subseteq Z(G)$ and
$(\beta,\gamma,\delta)\in\{(\m{yx},\m{xy},\m{xy})$, $(\m{xy},\m{yx},\m{yx})$,
$(\m{yx},\m{yx},\m{yx})\}$,

\item[(iv)] $G^4=1$ and $(\beta,\gamma,\delta) = (\m{x^{-1}y},\m{xy},\m{xy})$,

\item[(v)] $G^2\subseteq Z(G)$, $G^4=1$ and
$(\beta,\gamma,\delta)\in\{(\m{yx^{-1}},\m{xy},\m{xy})$,
$(\m{yx^{-1}},\m{yx},\m{yx})$, $(\m{x^{-1}y},\m{yx},\m{yx})\}$.
\end{enumerate}
\end{lemma}
\begin{proof}
We must check for which admissible values of $\alpha$, $\beta$, $\gamma$ and
$\delta$ the identities $(\ref{Eq:LA1})$--$(\ref{Eq:LA4})$ hold.

The identities $(\ref{Eq:LA1})$, $(\ref{Eq:LA2})$ always hold.

Denote by $I(\delta)$ the identity $(\ref{Eq:LA3})$, i.e.,
$\delta(x,\gamma(x,y))=\delta(x,x)y$. First observe that if $\delta$ inverts
its second argument then $I(\delta)$ does not have a solution; otherwise
$\gamma$ would have to invert its second argument, too, and that is not allowed
by Assumption \ref{As:4}. Now, $I(\m{xy})$ is $x\gamma(x,y)=xxy$, which has a
unique solution $\gamma=\m{xy}$; $I(\m{x^{-1}y})$ is $x^{-1}\gamma(x,y)=y$ with
solution $\gamma=\m{xy}$; $I(\m{yx})$ is $\gamma(x,y)x=x^2y$ with solution
$\gamma=\m{yx}$, but only if $G^2\subseteq Z(G)$; and $I(\m{yx^{-1}})$ is
$\gamma(x,y)x^{-1}=y$ with solution $\gamma=\m{yx}$.

We now consider the pairs $(\gamma,\delta)=(\m{xy},\m{xy})$,
$(\m{xy},\m{x^{-1}y})$, $(\m{yx},\m{yx})$, $(\m{yx},\m{yx^{-1}})$, and test for
which values of $\beta$ they satisfy $(\ref{Eq:LA4})$. It turns out that every
such triple $(\beta,\gamma,\delta)$ satisfies $(\ref{Eq:LA4})$, but additional
assumptions on $G$ are sometimes needed. Here is the calculation:

Assume $(\gamma,\delta)=(\m{xy},\m{xy})$, and denote by $J(\beta)$ the identity
$(\ref{Eq:LA4})$, i.e., $xxy=\beta(xx,y)$. Then $J(\m{xy})$ holds, $J(\m{yx})$
holds if $G^2\subseteq Z(G)$, $J(\m{yx^{-1}})$ is $x^2y=yx^{-2}$, which holds
if and only if $G^4=1$ and $G^2\subseteq Z(G)$, and $J(\m{x^{-1}y})$ holds when
$G^4=1$.

Assume $(\gamma,\delta)=(\m{xy},\m{x^{-1}y})$ and denote by $J(\beta)$ the
identity $(\ref{Eq:LA4})$, i.e., $y=\beta(1,y)$. Since $\beta$ never inverts
its second argument, $J(\beta)$ always holds.

Assume $(\gamma,\delta)=(\m{yx},\m{yx})$, $G^2\subseteq Z(G)$, and denote by
$J(\beta)$ the identity $(\ref{Eq:LA4})$, i.e., $yx^2=\beta(x^2,y)$. Then
$J(\m{xy})$ holds, $J(\m{yx})$ holds, $J(\m{yx^{-1}})$ holds if $G^4=1$, and
$J(\m{x^{-1}y})$ holds if $G^4=1$.

Assume $(\gamma,\delta)=(\m{yx},\m{yx^{-1}})$ and denote by $J(\beta)$ the
identity $(\ref{Eq:LA4})$, i.e., $y=\beta(1,y)$. Using our usual trick with the
second argument, we see that $J(\beta)$ always holds.
\end{proof}

\begin{theorem}\label{Th:Bol}
Suppose that $G$ is a nonabelian group and that the Assumptions
$\ref{As:1}$--$\ref{As:4}$ are satisfied. Then
$(M,*)=G(\alpha,\beta,\gamma,\delta)$ is a group if and only if $(M,*)$ is
equal to
\begin{equation}\label{Eq:XLB1}
    G(\m{xy},\m{xy},\m{xy},\m{xy}).
\end{equation}

The loop $(M,*)$ is a nonassociative Moufang loop if and only if it is equal to
\begin{equation}\label{Eq:XLB2}
    G(\m{xy},\m{x^{-1}y},\m{yx},\m{yx^{-1}}),
\end{equation}
and then it is isomorphic to the Chein loop $M(G,2)$.

The loop $(M,*)$ is a left Bol loop that is not Moufang if and only if: either
$G^2\subseteq Z(G)$ and $(M,*)$ is among
\begin{eqnarray}
    &&G(\m{xy},\m{x^{-1}y},\m{xy},\m{x^{-1}y}),\label{Eq:XLB3}\\
    &&G(\m{xy},\m{xy},\m{yx},\m{yx});\label{Eq:XLB4}
\end{eqnarray}
or $G^2\subseteq Z(G)$, $G^4=1$ and $(M,*)$ is among
\begin{eqnarray}
    &&G(\m{xy},\m{xy},\m{xy},\m{x^{-1}y}),\label{Eq:XLB5}\\
    &&G(\m{xy},\m{xy},\m{yx},\m{yx^{-1}}),\label{Eq:XLB6}\\
    &&G(\m{xy},\m{x^{-1}y},\m{xy},\m{xy}),\label{Eq:XLB7}\\
    &&G(\m{xy},\m{x^{-1}y},\m{yx},\m{yx}).\label{Eq:XLB8}
\end{eqnarray}
\end{theorem}
\begin{proof}
According to Lemma \ref{Lm:LA}, there are $16$ left alternative loops in the
nonabelian case. We now use the left Bol identity $(\ref{Eq:LB3})$ to eliminate
$8$ of them.

When $\gamma=\m{xy}$, $(\ref{Eq:LB3})$ becomes $\beta(x,yxz)=\beta(x,yx)z$, and
it is not satisfied when $\beta=\m{yx}$, $\beta=\m{yx^{-1}}$. When
$\gamma=\m{yx}$, $(\ref{Eq:LB3})$ becomes $\beta(x,xzy)=z\beta(x,xy)$, and it
is not satisfied when $\beta=\m{yx}$, $\beta=\m{yx^{-1}}$.

We claim that the remaining $8$ loops $G(\alpha,\beta,\gamma,\delta)$ are all
left Bol. It is clear for $(\ref{Eq:XLB1})$. The loop $(\ref{Eq:XLB2})$ is
isomorphic to $M(G,2)$ via Lemma \ref{Lm:Iso}. Since $G$ is nonabelian,
$(\ref{Eq:XLB2})$ is nonassociative \cite{Chein}.

A straightforward verification of identities $(\ref{Eq:LB1})$--$(\ref{Eq:LB8})$
shows that the loops $(\ref{Eq:XLB3})$--$(\ref{Eq:XLB8})$ are left Bol,
however, additional assumptions on $G$ are needed. We will indicate when the
assumptions are needed, and leave the verification of
$(\ref{Eq:LB1})$--$(\ref{Eq:LB8})$ to the reader.

First of all, $G^2\subseteq Z(G)$ is needed for $(\ref{Eq:XLB4})$ already in
Lemma \ref{Lm:LA}, and $G^2\subseteq Z(G)$, $G^4=1$ is needed for
$(\ref{Eq:XLB8})$ by the same Lemma. The condition $G^2\subseteq Z(G)$ is
required for $(\ref{Eq:XLB3})$ in $(\ref{Eq:LB5})$. The conditions
$G^2\subseteq Z(G)$, $G^4=1$ are required for $(\ref{Eq:XLB5})$ and
$(\ref{Eq:XLB6})$ in $(\ref{Eq:LB5})$, and for $(\ref{Eq:XLB7})$ in
$(\ref{Eq:LB6})$.

The loops $(\ref{Eq:XLB3})$--$(\ref{Eq:XLB8})$ are not Moufang, by Theorem
\ref{Th:Moufang}.
\end{proof}

\begin{remark} There are nonabelian groups satisfying $G^2\subseteq Z(G)$ and
$G^4=1$. For instance the $8$-element dihedral group $D_8$ and the $8$-element
quaternion group $Q_8$ have this property. So do the groups $(D_8)^a\times
(Q_8)^b\times (C_4)^c \times (C_2)^d$ with $a+b>0$.
\end{remark}

\section{Questions}

\noindent Our methods certainly do not yield all Bol loops with a subgroup of
index $2$. This is apparent already from Remark \ref{Rm:Middle}. Is there some
way of determining all such Bol loops?

Is it true that the $6$ constructions $(\ref{Eq:XLB3})$--$(\ref{Eq:XLB8})$ of
Theorem \ref{Th:Bol} produce $6$ pairwise nonisomorphic loops when $G$ is
fixed? Note that $\oo{x}\in\oo{G}$ is of order $2$ if and only if
$\delta(x,x)=1$. Hence there are $|G_2|+|G|$ elements of exponent $2$ in
$(\ref{Eq:XLB3})$, $(\ref{Eq:XLB5})$, $(\ref{Eq:XLB6})$, and $2|G_2|<|G|+|G_2|$
elements of exponent $2$ in the remaining three loops. The following lemma
allows us to distinguish more loops:

\begin{lemma} Let $(\ref{Eq:XLB4})$, $(\ref{Eq:XLB7})$, $(\ref{Eq:XLB8})$ be as
in Theorem $\ref{Th:Bol}$. Then $(\ref{Eq:XLB7})\not\cong(\ref{Eq:XLB4})
\not\cong (\ref{Eq:XLB8})$.
\end{lemma}
\begin{proof}
Let $(M,*)=G(\alpha,\beta,\gamma,\delta)$. Let us count the cardinality of
$Z(M)=\{x\in M;\;x*y=y*x$  for every $y\in M\}$ for the three loops in
question.

Assume $x\in G$. When $y\in G$, we have $x*y=y*x$ if and only if $x\in Z(G)$.
When $\oo{y}\in \oo{G}$, we have $x*\oo{y}=\oo{y}*x$ if and only if
$\beta(x,y)=\gamma(y,x)$. This is always true for $(\ref{Eq:XLB4})$.

Assume $\oo{x}\in\oo{G}$.  When $\oo{y}\in\oo{G}$, we have
$\oo{x}*\oo{y}=\oo{y}*\oo{x}$ if and only if $\delta(x,y)=\delta(y,x)$. This is
true (in all three cases) if and only if $x\in Z(G)$. When $x\in Z(G)$ and
$y\in G$, we have $\oo{x}*y=y*\oo{x}$ if and only if $\gamma(x,y)=\beta(y,x)$.
This reduced to $yx=yx$ (always true) for $(\ref{Eq:XLB4})$, to $xy=y^{-1}x$
(always false, because $xy=xy^{-1}$ holds for every $y$ if and only if $G_2=G$)
for $(\ref{Eq:XLB7})$, and to $yx=y^{-1}x$ (always false) for
$(\ref{Eq:XLB8})$.

Altogether, $|Z(M)| = 2|Z(G)|$ for $(\ref{Eq:XLB4})$, while $|Z(M)| \le |Z(G)|$
for the other two loops.
\end{proof}

Even if the $6$ loops are pairwise nonisomorphic for a given $G$, it is
possible that $G(\alpha,\beta,\gamma,\delta)\cong
H(\alpha',\beta',\gamma',\delta')$ holds for some nonisomorphic groups $G$, $H$
to which Theorem \ref{Th:Bol} is applied. Can you determine these
``exceptional'' isomorphisms? Are they exceptional or common?

\bibliographystyle{plain}

\end{document}